\documentclass[11pt]{article}
\usepackage{amsmath,amssymb}
\begin{document}
\title{\bf States on sharply dominating effect algebras}
\author{{Zdenka Rie\v{c}anov\'a}
\thanks {Department of Mathematics, Faculty of Electrical
Engineering and Information Technology STU, Ilkovi\v{c}ova~3,
SK-812~19~Bratislava, Slovakia. E-mail:
zdenka.riecanova@stuba.sk}, {Wu Junde}\thanks {Corresponding
Author Wu Junde: Department of Mathematics, Zhejiang University,
Hangzhou 310027, People's Republic of China. E-mail:
wjd@zju.edu.cn}}
\date{}
\maketitle

\begin{tabular*}{13cm}{r}
\hline
\end{tabular*}

{\small\it We prove that Archimedean sharply dominating atomic
lattice effect algebras can be characterized by property called
``basic decomposition of elements". As an application we prove the
``state smearing theorem" for these effect algebras.}

\begin{tabular*}{13cm}{r}
\hline
\end{tabular*}



\vskip 0.1 in

{\bf 1. Introduction}

\vskip 0.1 in

In recent years quantum effects and fuzzy events are studied
within a general algebraic framework called effect algebras
introduced by Foulis and Bennett [2], or equivalent in some sense
$D$-posets [8], or $D$-algebra [4]. Lattice ordered effect
algebras generalize orthomodular lattices (which may include
noncompatible pairs of elements) and $MV$-algebras, originally
constructed to give an algebraic structure to the infinite-valued
Lukasiewicz propositional logics [1]. Thus effect algebras have
their importance in the investigation of the phenomenon of
uncertainty.

All these generalizations of Boolean algebras are investigated as
carriers of states and probability measures. Nevertheless, there
are even (finite) effect algebras admitting no states (see [14],
[3]). In [7] has been shown that the set of sharp elements $S(E)$
of a lattice effect algebra $E$ is an orthomodular lattice which
is a sublattice and a sub-effect algebra of $E$. In this sense,
lattice effect algebras are smeared orthomodular lattices of their
sharp elements. Simultaneously, in some cases of complete atomic
effect algebras, states existing on the set $S(E)$ of sharp
elements can be smeared onto whole effect algebra $E$, [16], [18].

Unfortunately, not every (lattice) effect algebra can be embedded
into a complete effect algebra (see [11]). Hence, we are
interested in finding a family of atomic lattice effect algebras
which elements have a ``basic decomposition" into a sum of a
unique sharp element and a unique $\oplus$-orthogonal set of
unsharp elements being multiplies of atoms. We have found that
this claim completely characterize the family of Archimedean
sharply dominating lattice effect algebras. The notion of sharply
dominating (S-dominating) effect algebras has been introduced by
S. Gudder (see [5], [6]). As an important application, we obtain
the ``state smearing theorem" for this family of effect algebras.

\vskip 0.1 in

{\bf 2. Effect algebras, basic notions and facts}

\vskip 0.1 in

Effect algebras have been introduced by Foulis and Bennett [2] as
an algebraic structure pvoviding an instrument for studying
quantum effects that may be unsharp.

{\bf Definition 2.1.} A structure  $(E; \oplus, 0, 1)$ is called
an {\it effect algebra} if 0, 1 are two distinguished elements and
$\oplus$ is a partially defined binary operation on $E$ which
satisfies the following conditions for any $a, b, c\in E$:

(Ei). $a\oplus b=b\oplus a$ if $a\oplus b$ is defined,

(Eii). $(a\oplus b)\oplus c=a\oplus (b\oplus c)$ if one side is
defined,

(Eiii). For each $a\in E$ there exists a unique $b\in E$ such that
$a\oplus b=1$ (we put $a^{'}=b$),

(Eiv). If $1\oplus a$ is defined then $a=0$.

We often denote the effect algebra $(E; \oplus, 0, 1)$ briefly by
$E$. In every effect algebra $E$ we can define the partial
operation $\ominus$ and the partial order $\leq$ by putting $a\leq
b$ and $b\ominus a=c$ iff $a\oplus c$ is defined and $a\oplus
c=b$.

Note that for $a, b\in E, a\oplus b$ is defined iff $b\leq a^{'}$.

If $E$ with the defined partial order is a lattice (complete
lattice), then $E$ is called a lattice (complete) effect algebra.

Recall that $Q\subseteq E$ is called a {\it sub-effect algebra} of
$E$ iff

(i) $1\in Q$,

(ii) if $a, b, c\in E$ with $a\oplus b=c$ and out of $a, b, c$ at
least two elements are in $Q$ then $a, b, c\in Q$.

In this case, $Q$ with inherited operation $\oplus$ is an effect
algebra in its own right.

Elements $a, b$ of a lattice effect algebra $(E; \oplus, 0, 1)$
are called {\it compatible} (written $a\leftrightarrow b$) iff
$a\vee b=a\oplus (b\ominus(a\wedge b))$. If for every $a, b\in E$
we have $a\leftrightarrow b$ then $E$ is an $MV$-effect algebra,
[9].

It is important to note that every orthomodular lattice $(L; \vee,
\wedge, ', 0, 1)$ can be itself organized into a lattice effect
algebra if we define: $a\oplus b=a\vee b$ iff $a\leq b^{'}$. Then
$(L; \oplus, 0, 1)$ is a lattice effect algebra.

{\bf Lemma 2.2.} Let $(E; \oplus, 0, 1)$ be a lattice effect
algebra and $x, y, z\in E$.

(i) If $x\leq y^{'}$ then $x\oplus y=(x\vee y)\oplus (x\wedge y)$.

(ii) If $x, y\leq z^{'}$ then $(x\vee y)\oplus z=(x\oplus z)\vee
(y\oplus z)$.

(iii) If $x\wedge y=0$ and for $m, n\in N$ the elements $mx, ny$
and $(mx)\oplus (ny)$ exist in $E$ then $(kx)\wedge (ly)=0$ and
$(kx)\vee (ly)=(kx)\oplus (ly)$ for all $k\in \{1, 2, ..., m\},
l\in \{1, 2, ..., n\}.$

(iv). If $Y\subseteq E$ with $\bigvee Y$ existing in $E$ and $x\in
E$ such that $x\leftrightarrow y$ for all $y\in Y$ then $x\wedge
\bigvee Y=\bigvee \{x\wedge y|y\in Y\}$ and $x\leftrightarrow
\bigvee Y$.

For the proof of (i) we refer to [10], the proof of (ii) is left
to the reader, for the proof of (iii) we refer to [17] and (iv) to
[7].

We say that a finite system $F=(a_k)_{k=1}^n$ of not necessarily
different elements of an effect algebra $(E; \oplus, 0, 1)$ is
$\oplus$-orthogonal if $a_1\oplus a_2\oplus a_3...\oplus a_n$
(written $\bigoplus_{k=1}^na_k$ or $\bigoplus F$) exists in $E$.
Here we define $a_1\oplus a_2\oplus a_3...\oplus a_n=(a_1\oplus
a_2\oplus a_3...\oplus a_{n-1})\oplus a_n$ suppose that
$\oplus_{k=1}^{n-1}a_k$ exists and $\oplus_{k=1}^{n-1}a_k\leq
a^{'}_n$. An arbitrary system $G=(a_{\kappa})_{\kappa\in H}$ of
not necessarily different elements of $E$ is called
$\oplus$-orthogonal if $\bigoplus K$ exists for every finite
$K\subseteq G$. We say that for a $\oplus$-orthogonal system
$G=(a_{\kappa})_{\kappa\in H}$ the element $\bigoplus G$ exists
iff $\bigvee\{\bigoplus K|K\subseteq G$ finite $\}$ exists in $E$
and then we put $\bigoplus G = \bigvee\{\bigoplus K|K\subseteq G$
finite $\}$ (we write $G_1\subseteq G$ iff there is $H_1\subseteq
H$ such that $G_1=(a_{\kappa})_{\kappa\in H_1}$).

An effect algebra $(E; \oplus, 0, 1)$ is called {\it Archimedean}
if for no nonzero element $e\in E$ the elements $ne=e\oplus
e\oplus...\oplus e$ (n times) exist for all $n\in N$. We can show
that every complete effect algebra is Archimedean [12].

For an element $x$ of an effect algebra $E$ we write
ord$(x)=\infty$ if $nx$ exists for every $n\in N$. We write
ord$(x)=n_x\in N$ if $n_x$ (called isotropic index of $x$) is the
greatest integer such that $n_xx$ exists in $E$. Clearly, in an
Archimedean effect algebra $n_x< \infty$ for every $x\in E$.

Recall that $x\in E$ is called a {\it sharp} element of an effect
algebra $E$ if $x\wedge x^{'}=0$. The set $S(E)=\{x\in E| x\wedge
x^{'}=0\}$ is called a set of all sharp elements of $E$ [5]. It
has been shown in [7] that in every lattice effect algebra $E$ the
set $S(E)$ is an orthomodular lattice, being a sub-effect algebra
and a sublattice of $E$. Moreover, $S(E)$ is a full sublattice of
$E$, which means that $S(E)$ inherits all suprema and infima of
subsets of $S(E)$, existing in $E$.

A nonzero element $a$ of an effect algebra $E$ is called an {\it
atom} if $0\leq b< a$ implies $b=0$. $E$ is called {\it atomic} if
for every nonzero element $x\in E$ there is an atom $a\in E$ such
that $a\leq x$.

{\bf Lemma 2.3.} Assume that $(E; \oplus, 0, 1)$ is an atomic
lattice effect algebra and $a\in E$ is an atom with ord$(a)=n_a\in
N$. Let $S(E)=\{x\in E|x\wedge x^{'}=0\}$. Then

(i) $(ka)\wedge (ka)^{'}\neq 0$ for all $k\in \{1, 2, ...,
n_a-1\}$,

(ii) $n_aa\in S(E)$ and $ka\notin S(E)$ for $k\neq n_a$.

(iii) if $x\in E$ with $a\leq x\leq ka$ then there is $r\in N$
such that $x=ra$,

(iv) If $a, b\in E$ are atoms and $k, l\in N$ are such that $k\neq
n_a$ and $ka=lb$ then $a=b$ and $k=l$.

(v) If $E$ is Archimedean then to every $x\in E, x\neq 0$ there
exist a set $\{a_{\alpha}|\alpha\in \Lambda\}$ of atoms and
positive integers $k_{\alpha}$ such that

$$x=\bigoplus \{k_{\alpha}a_{\alpha}|\alpha\in \Lambda\}=\bigvee
\{k_{\alpha}a_{\alpha}|\alpha\in \Lambda\},$$ under which $x\in
S(E)$ iff $k_{\alpha}=n_{a_{\alpha}}=$ ord$(a_{\alpha})$ for all
$\alpha\in\Lambda$.

For the proof we refer the readers to [17] Theorem 2.4 and [16],
Theorem 3.3.

{\bf Theorem 2.4.} Let $(E; \oplus, 0, 1)$ be a lattice effect
algebra. Let for atoms $a, b\in E, ka\leq lb$. Then

(i) If ord$(b)=\infty$ or $l<$ord$(b)$, then $a=b$.

(ii) If $l=$ord$(b)=n_b<\infty$, then either $a=b$ or
$a\nleftrightarrow b$ ($a, b$ are noncompatible) and $n_aa\leq
n_bb$.

{\bf Proof.} (i) Assume that $a\neq b$.  By assumptions we have
$a\leq ka\leq lb\leq b^{'}$, which, by Lemma 2.2 implies that
$a\oplus b=(a\vee b)\oplus (a\wedge b)=a\vee b\leq b^{'}$ and
hence $a\oplus 2b=(a\vee b)\oplus b=(a\oplus b)\vee 2b=a\vee
2b\leq b^{'}$. By induction it follows that $a\oplus lb=a\vee
lb=lb$ which gives $a=0$, a contradiction. Hence $a=b$.

(ii) Assume that $a\neq b$ and $a\leftrightarrow b$. Then $a\vee
b=a\oplus (b\ominus (a\wedge b))=a\oplus b\leq b^{'}$, which gives
that $a\oplus 2b=a\vee 2b\leq b^{'}$. By induction $a\oplus
lb=a\vee lb=lb$ and hence $a=0$, a contradiction. Hence
$a\nleftrightarrow b$. Moreover, $a\leq n_bb$ implies the
existence of $a\oplus (n_bb)^{'}=(a\vee (n_bb)^{'})\oplus (a\wedge
(n_bb)^{'})=a\vee (n_bb)^{'}\leq a^{'}$ and by induction,
$(n_aa)\oplus (n_bb)^{'}$ exists, which gives $n_aa\leq n_bb$.

Note that if effect algebra $E$ is not lattice ordered, then
statements of Lemma 2.3 and Theorem 2.4 fails to be true in
general.

{\bf Example 2.5.} Let $E=\{0, a, b, a\oplus b, 2a=2b, 1=2a\oplus
b=3b\}$. Then ord$(a)=2$, ord$(b)=3$ and $S(E)=\{0, 1\}$. Hence
$n_aa=2a\notin S(E)$. Moreover, $2b=2a=n_aa$ in spite of that
$a\neq b$ and $a\leftrightarrow b$ because $a\oplus b$ exists.

{\bf Theorem 2.6.} Let  $(E; \oplus, 0, 1)$ be a lattice effect
algebra. For any two sets $\{a_{\alpha}|\alpha\in\Lambda\}$,
$\{b_{\beta}|\beta\in \cal B\}$ of atoms of $E$ and positive
integers $k_{\alpha}\neq$ ord$(a_{\alpha})$, $l_{\beta}\neq$
ord$(b_{\beta})$, the following conditions are equivalent:

(i)
$\bigoplus\{k_{\alpha}a_{\alpha}|\alpha\in\Lambda\}=\bigoplus\{l_{\beta}b_{\beta}|\beta\in
\cal B\}$.

(ii) For every $\alpha\in\Lambda$ there exists $\beta\in\cal B$
such that $a_{\alpha}=b_{\beta}$ and $k_{\alpha}=l_{\beta}$.

{\bf Proof.} Let $\alpha_0\in\Lambda$. Then
$k_{\alpha_0}a_{\alpha_0}\leq a^{'}_{\alpha_0}$ because
$k_{\alpha_0} < $ ord$(a_{\alpha_0})$. Moreover, for every
$\alpha\neq \alpha_0, \alpha\in \Lambda$ there exists
$k_{\alpha}a_{\alpha}\oplus k_{\alpha_0}a_{\alpha_0}$ which gives
$k_{\alpha}a_{\alpha}\leq (k_{\alpha_0}a_{\alpha_0})^{'}\leq
a^{'}_{\alpha_0}$. It follows that for every $\beta\in\cal B$ we
have $$b_{\beta}\leq
l_{\beta}b_{\beta}\leq\bigoplus\{k_{\alpha}a_{\alpha}|\alpha\in\Lambda\}=\bigvee
\{k_{\alpha}a_{\alpha}|\alpha\in\Lambda\}\leq a^{'}_{\alpha_0}.$$
It follows that $b_{\beta}\leftrightarrow a_{\alpha_0}$ and by
[13] we obtain that $k_{\alpha_0}a_{\alpha_0}\leftrightarrow
l_{\beta}b_{\beta}$ for all $\beta\in\cal B$. By Lemma 2.2, (iv),

$k_{\alpha_0}a_{\alpha_0}=k_{\alpha_0}a_{\alpha_0}\wedge (\bigvee
\{l_{\beta}b_{\beta}|\beta\in{\cal
B}\})=\bigvee\{k_{\alpha_0}a_{\alpha_0}\wedge
l_{\beta}b_{\beta}|\beta\in\cal B\}$.

Furthermore, by Lemma 2.2, (iii), $k_{\alpha_0}a_{\alpha_0}\wedge
l_{\beta}b_{\beta}=0$ for all $b_{\beta}\neq a_{\alpha_0}$ and
hence there exists $\beta_0\in\cal B$ such that
$a_{\alpha_0}=b_{\beta_0}$ and $k_{\alpha_0}a_{\alpha_0}=
l_{\beta_0}a_{{\alpha_0}}$ which gives $k_{\alpha_0}=l_{\beta_0}$.

\vskip 0.1 in

{\bf 3. Sharply dominating lattice effect algebras}

\vskip 0.1 in

A special types of effect algebras called sharply dominating and
$S$-dominating has been introduced by S. Gudder in [5] and [6].
Important example is a standard Hilbert spaces effect algebra
$\cal E(H)$ of positive linear operators on a complex Hilbert
space $\cal H$ which are dominated by identity operator $I$. See
[5] and [6].

{\bf Definition 3.1} ([5], [6]). An effect algebra $(E, \oplus, 0,
1)$ is called sharply dominating if for every $a\in E$ there
exists a smallest sharp element $\hat{a}$ such that $a\leq
\hat{a}$. That is $\hat{a}\in S(E)$ and if $b\in S(E)$ satisfies
$a\leq b$ then $\hat{a}\leq b$. A {\it sharply dominating effect
algebra} $E$ is called $S$-{\it dominating} if $a\wedge p$ exists
for every $a\in E$ and $p\in S(E)$.

In next we will use that $S(E)$ is a sub-lattice and a sub-effect
algebra of $E$, [7].

Clearly a lattice effect algebra is $S$-dominating iff $E$ is
sharply dominating.

{\bf Lemma 3.2.} Let $(E; \oplus, 0, 1)$ be a lattice effect
algebra. The following conditions are equivalent:

(i). $E$ is sharply dominating.

(ii). For every $x\in E$ there exists $\tilde{x}\in S(E)$ such
that $\tilde{x}\leq x$ and if $u\in S(E)$

satisfies $u\leq x$ then $u\leq \tilde{x}$.

(iii). For every $x\in E$ there exists the unique element $v_x\in
S(E)$ satisfying $v_x\leq x$

and if $v\in S(E)$ satisfies $v\leq x\ominus v_x$ then $v=0$.

{\bf Proof.} (i) $\Leftrightarrow$ (ii). This is easy to verify by
deMorgan's laws, since $(E, \leq, 0, 1)$ is a DM-poset.

(ii) $\Longrightarrow$ (iii). Let $\tilde{x}$ be a greatest sharp
element under $x$. If $v\in S(E)$ satisfies $v\leq x\ominus
\tilde{x}$ then $v\oplus \tilde{x}\leq x$ and because $v\oplus
\tilde{x}\in S(E)$ we obtain that $v\oplus \tilde{x}\leq
\tilde{x}$ which gives $v=0$. Assume that $u\in S(E)$, $u\leq x$
and if $v\in S(E)$ satisfies $v\leq x\ominus u$ then $v=0$. Then
$u\leq \tilde{x}\leq x$ and hence $\tilde{x}\ominus u\leq x\ominus
u$, which implies that $\tilde{x}\ominus u=0$, because by [7]
$\tilde{x}\ominus u\in S(E)$. This proves that $u=\tilde{x}$,
hence $v_x=\tilde{x}$ is the unique element satisfying (iii).

(iii) $\Longrightarrow$ (ii). Assume that $x\in E, x\neq 0$ and
$v_x$ satisfies (iii). Then $v_x\in S(E)$ and $v_x\leq x$.
Further, if $u\in S(E)$ satisfies $u\leq x$ then by [7] $v_x\vee
u\in S(E)$ and $v_x\vee u\leq x$. Moreover, if $v\leq x\ominus
(v_x\vee u)\leq x\ominus v_x$ then $v=0$, which by (iii) gives
that $v_x\vee u=v_x$ and hence $u\leq v_x$ which proves that
$v_x=x$ satisfies (ii).

{\bf Remark 3.3.} There are atomic lattice effect algebras (even
MV-effect algebras) which are sharply dominating but non
Archimedean. The simplest example is the Chang MV-effect algebra
$E=\{0, a, 2a, 3a,..., (3a)^{'}, (2a)^{'}, a^{'}, 1\}$, which has
the unique atom $a$ with $ord(a)=\infty$. The set of sharp
elements is $S(E)=\{0, 1\}$, hence $E$ is sharply dominating.

The next theorem give a characterization of all atomic lattice
effect algebras which are sharply dominating and Archimedean. In
what follows set

$M(E)=\{x\in E|$ if $v\in S(E)$ satisfies $v\leq x$ then $v=0\}.$

{\bf Theorem 3.4.} Let $(E; \oplus, 0, 1)$ be an atomic lattice
effect algebra. The following conditions are equivalent:

(i). $E$ is Archimedean and sharply dominating.

(ii). For every $x\in E, x\neq 0$ there exists the unique $v_x\in
S(E)$, unique set of atoms $\{a_{\alpha}|\alpha\in\Lambda\}$ and
unique positive integers $k_{\alpha}\neq ord(a_{\alpha})$ such
that
$$x=v_x\oplus (\bigoplus\{k_{\alpha}a_{\alpha}|\alpha\in\Lambda\}).$$

{\bf Proof.} (i) $\Longrightarrow$ (ii). Let $x\in E, x\neq 0$. By
Lemma 3.2 there exists the unique $v_x\in S(E)$ such that
$x=v_x\oplus (x\ominus v_x)$ and $x\ominus v_x\in M(E)$. By [16],
Theorem 3.3 there exist a set $\{a_{\alpha}|\alpha\in\Lambda\}$ of
atoms and positive integers  $k_{\alpha}, \alpha\in\Lambda$ such
that $x\ominus
v_x=\bigoplus\{k_{\alpha}a_{\alpha}|\alpha\in\Lambda\}=\bigvee\{k_{\alpha}a_{\alpha}|\alpha\in\Lambda\}$.
Because $n_{a_{\alpha}}a_{\alpha}\in S(E)$ and $x\ominus v_x\in
M(E)$, we obtain that $k_{\alpha}\neq n_{a_{\alpha}}$. By Theorem
2.6 the set $\{a_{\alpha}|\alpha\in\Lambda\}$ and positive
integers $k_{\alpha}, \alpha\in\Lambda$ are unique.

(ii) $\Longrightarrow$ (i). Let $x\in E, x\neq 0$ and let $v_x$
and the set of atoms $\{a_{\alpha}|\alpha\in\Lambda\}$ satisfy
(ii). Clearly if $x\in S(E)$ then $v_x=x=\hat{x}$ is the smallest
sharp element dominating $x$. Assume that $x\notin S(E)$. Then
$x\ominus v_x\in M(E)$, because otherwise there exists $v\in S(E),
v\neq 0$ with $v\leq x\ominus v_x$. Further, $v\neq x\ominus v_x$
since $x\notin S(E)$. It follows that $x=v_x\oplus v\oplus
(x\ominus(v_x\oplus v))$ and by (ii) there is $u\in S(E)$ and the
set of atoms $\{b_{\beta}|\beta\in \cal B\}$ such that $x\ominus
(v_x\oplus v)=u\oplus (\bigoplus \{b_{\beta}|\beta\in \cal B\})$,
which contradicts to (ii). By Lemma 3.2 we obtain that $E$ is
sharply dominating.

Further, let us show that for every atom $a\in E$ we have
$ord(a)<\infty$. Assume to the contrary that there exists an atom
$a\in E$ with $ord(a)=\infty$. Let $\omega=\hat {a}\in S(E)$ is
the smallest sharp element dominating $a$. Then $\omega=a\oplus
(\omega\ominus a)$ and $\omega\ominus a\in M(E)$, because
otherwise there is $v\neq 0, v\in S(E)$ which satisfies $v\leq
\omega\ominus a$ and hence $a\leq \omega\ominus v\in S(E)$ and
$\omega\leq \omega\ominus v$, a contradiction. By (ii) there exist
a unique set of atoms $\{a_{\alpha}|\alpha\in\Lambda\}$ and
positive integers $k_{\alpha}, k_{\alpha}\neq ord(a_{\alpha}),
\alpha\in\Lambda$ such that $\omega\ominus a=
\bigoplus\{k_{\alpha}a_{\alpha}|\alpha\in\Lambda\}=\bigvee\{k_{\alpha}a_{\alpha}|\alpha\in\Lambda\}$,
hence $\omega =
a\oplus\bigvee\{k_{\alpha}a_{\alpha}|\alpha\in\Lambda\}=\bigvee\{(k_{\alpha}a_{\alpha})\oplus
a|\alpha\in\Lambda\}$.

If there exists $\alpha_0\in\Lambda$ such that $a=a_{\alpha_0}$
then for every $\alpha\in\Lambda, \alpha\neq \alpha_0$ we have
$k_{\alpha}a_{\alpha}\oplus a_{\alpha_0}=k_{\alpha}a_{\alpha}\vee
a_{\alpha_0}\leq a^{'}_{\alpha_0}$, since $a_{\alpha}\wedge
a_{\alpha_0}=0$ (see Lemma 2.2, (iii)). Because
$ord(a_{\alpha_0})=ord(a)=\infty$, we have
$k_{\alpha_0}a_{\alpha_0}\oplus a=k_{\alpha_0}a_{\alpha_0}\oplus
a_{\alpha_0}=(k_{\alpha_0}+1)a_{\alpha_0}<a^{'}_{\alpha_0}$. It
follows that $a\leq\omega\leq a^{'}_{\alpha_0}=a^{'}$, hence
$a\leq \omega^{'}$, a contradiction. If $a\notin
\{a_{\alpha}|\alpha\in\Lambda\}$, then by Lemma 2.2, (iv), we have
$a\wedge (\bigvee \{k_{\alpha}a_{\alpha}|\alpha\in\Lambda\}=0$ and
hence, by Lemma 2.2, (i), we have $a\leq \omega\leq a\vee (\bigvee
\{k_{\alpha}a_{\alpha}|\alpha\in\Lambda\})\leq a^{'}$, a
contradiction.

We have proved that $ord(a)=n_{\alpha}< \infty$, for every atom
$a\in E$. It follows that for every $x\in E, x\neq 0$ we have
$ord(x)=n_x< \infty$. Otherwise, since $E$ is atomic, there exists
an atom $a\in E, a\leq x$ and hence for every positive integer
$n$, elements $na\leq nx$ exist, hence $ord(a)=\infty$, a
contradiction. This proves that $E$ is Archimedean.

{\bf Theorem 3.5.} Let $(E; \oplus, 0, 1)$ be a lattice effect
algebra. For every atom $a\in E$ with $ord(a) < \infty$, $n_aa$ is
the smallest sharp element dominating atom $a$.

{\bf Proof.} Let $a\in E$ is an atom and $\omega\in S(E)$ with
$a\leq\omega$. Then $a\oplus \omega^{'}$ is defined and because
$a\wedge\omega^{'}\leq \omega\wedge \omega^{'}=0$, we obtain that
$a\oplus \omega^{'}=a\vee \omega^{'}\leq a^{'}$, by Lemma 2.2,
(i). It follows that there exist $2a\oplus \omega^{'}=2a\vee
\omega^{'}\leq a^{'}$, by Lemma 2.2, (iii). Hence $3a\oplus\omega$
exists. By induction $n_aa\oplus \omega^{'}$ exists and hence
$n_aa\leq \omega$. This proves that $n_aa$ is the smallest sharp
element dominating $a$.

{\bf Remarks 3.6.} Statements of Theorems 3.4 and 3.5 fails to be
true if $E$ is not lattice ordered.

{\bf Example 3.7.} Consider the effect algebra $(E; \oplus, 0,
1)$, where $E=\{0, a, b, a\oplus b, 2a=2b, 3b=2a\oplus b=1\}$.

Since $E$ is finite and $S(E)=\{0, 1\}$, $E$ is Archimedean and
sharply dominating. In spite of $ord(a)=2$, we have $2a\notin
S(E)$ because $(2a)^{'}=b$ and hence $(2a)\wedge (2a)^{'}=b$. In
spite of $ord(a)=2<ord(b)=3$, we have $2a=2b$. Clearly $E$ is not
lattice ordered because $a\vee b$ does not exists in $E$.

\vskip 0.1 in

{\bf 4. Smearing of states on sharply dominating effect algebras}

\vskip 0.1 in

In [7] has been shown that the subset $S(E)=\{x\in E|x\wedge
x^{'}=0\}$ of a lattice effect algebra $E$ is an orthomodular
lattice that is sublattice and a sub-effect algebra of $E$. We are
going to show, using Theorem 3.4 on basic decomposition of
elements, that the state on sharply dominating Archimedean atomic
lattice effect algebra $E$ exists if there exists an
$(o)$-continuous state on sharp elements of $E$. Note that an
example of an effect algebra admitting no states has been
presented in [14].

A net $(a_{\alpha})_{\alpha\in\Lambda}$ of elements of a poset
$(P, \leq )$ {\it order converges} to a point $a\in P$ if there
are nets $(u_{\alpha})_{\alpha\in\Lambda}$ and
$(v_{\alpha})_{\alpha\in\Lambda}$ of elements of $P$ such that
$$a\uparrow u_{\alpha}\leq a_{\alpha}\leq v_{\alpha}\downarrow
a.$$ We write $a_{\alpha}\xrightarrow{(o)} a$ in $P$ (or briefly
$a_{\alpha}\xrightarrow{(o)} a$). Here $u_{\alpha}\uparrow a$
means that $u_{\alpha}\leq u_{\beta}$ for all $\alpha\leq \beta$
and $a=\bigvee \{u_{\alpha}|\alpha\in \Lambda\}$. The meaning of
$v_{\alpha}\downarrow a$ is dual.

Recall that a map $\omega: E\rightarrow [0, 1]$ is called a
(finite additive) {\it state} on an effect algebra $(E; \oplus, 0,
1)$ if $\omega (1)=1$ and $x\leq y^{'}\Rightarrow \omega (x\oplus
y)=\omega (x)+ \omega (y)$. A state is {\it faithful} if $\omega
(x)=0\Rightarrow x=0$. A state $\omega$ is called $(o)$-{\it
continuous (order-continuous)} if $x_{\alpha}\xrightarrow{(o)}
a\Rightarrow \omega(x_{\alpha})\rightarrow \omega (x)$ for every
net $(x_{\alpha})_{\alpha\in\Lambda}$ of elements of $E$.

{\bf Lemma 4.1}. A state $\omega$ on an effect algebra $E$ is
$(o)$-continuous iff $x_{\alpha}\downarrow 0\Rightarrow \omega
(x_{\alpha})\downarrow 0$ for $x_{\alpha}\in E$.

For a proof we refer the reader to ([15], Lemma 4.4]).

Finally, recall that a map $\omega: L\rightarrow [0, 1]$ is a
state on an orthomodular lattice $(L; \vee, \wedge, {'}, 0, 1)$
iff $\omega (x\vee y)=\omega (x)+\omega (y)$ for all $x\leq y^{'},
x, y\in L$. Since for lattice effect algebra $(L; \oplus, 0, 1)$
derived from the orthomodular lattice $L$ we have $x\oplus y=x\vee
y$ iff $x\leq y^{'}$, we conclude that $\omega$ is also a state on
the effect algebra $L$, and conversely.

{\bf Theorem 4.1}. Let $(E; \oplus, 0, 1)$ be a sharply dominating
Archimedean atomic lattice effect algebra. Let $x\in E$ and
$\{a_{\alpha}|\alpha\in \Lambda\}$ be a set of atoms of $E$ such
that $x=\bigoplus\{k_{\alpha}a_{\alpha}|\alpha\in\Lambda\}$. Then

(i). $\bigoplus\{k_{\alpha}a_{\alpha}|k_{\alpha}=ord(a_{\alpha}),
\alpha\in\Lambda\}$ and
$\bigoplus\{k_{\alpha}a_{\alpha}|k_{\alpha}\neq ord(a_{\alpha}),
\alpha\in\Lambda\}$ exist in $E$.

(ii).
$v_x=\bigoplus\{k_{\alpha}a_{\alpha}|k_{\alpha}=ord(a_{\alpha}),
\alpha\in\Lambda\}\in S(E)$ and $x=v_x\oplus(
\bigoplus\{k_{\alpha}a_{\alpha}|k_{\alpha}\neq ord(a_{\alpha}),
\alpha\in\Lambda\})$ is a basic decomposition of $x$.

{\bf Proof.} By Lemma 2.2, (iii) and the definition of $\bigoplus
D$ for a $\oplus$-orthogonal system $D$ we have
$\bigoplus\{k_{\alpha}a_{\alpha}|\alpha\in\Lambda\}=\bigvee\{\bigoplus\{k_{\alpha}a_{\alpha}|\alpha\in
K\}|K\subseteq \Lambda$ finite
$\}=\bigvee\{\bigvee\{k_{\alpha}a_{\alpha}|\alpha\in
K\}|K\subseteq \Lambda$ finite
$\}=\bigvee\{k_{\alpha}a_{\alpha}|\alpha\in\Lambda\}$. By Lemma
3.2 there exist unique elements $v_x\in S(E)$ and $x\ominus v_x\in
M(E)$ such that $x=v_x\oplus(x\ominus v_x)$. Since $x\ominus
v_x\leq v^{'}_x$ and $v_x\in S(E)$ we obtain that $v_x\wedge
(x\ominus v_x)=0$ and hence $v_x\oplus (x\ominus v_x)=v_x\vee
(x\ominus v_x)$ by Lemma 2.2, (i).

Let $\alpha_0\in\Lambda$ with $k_{\alpha_0}\neq
ord(a_{\alpha_0})$. Then $k_{\alpha_0}a_{\alpha_0}\leq
a^{'}_{\alpha_0}$ and for every $\alpha\in\Lambda, \alpha\neq
\alpha_0$ we have $k_{\alpha}a_{\alpha}\leq
(k_{\alpha_0}a_{\alpha_0})^{'}\leq a^{'}_{\alpha_0}$ since
$(k_{\alpha_0}a_{\alpha_0})\oplus k_{\alpha}a_{\alpha}$ is
defined. It follows that $v_x\oplus(x\ominus
v_x)=x=\bigvee\{k_{\alpha}a_{\alpha}|\alpha\in\Lambda\}\leq
a^{'}_{\alpha_0}$ which implies that $v_x$ and $x\ominus v_x$ are
compatible with $a_{\alpha_0}$ and hence also with
$k_{\alpha_0}a_{\alpha_0}$ (see [13]). Moreover,
$(k_{\alpha_0}a_{\alpha_0})\wedge v_x=0$ because otherwise by
Theorem 2.4 resp. its Corollary we have $a_{\alpha_0}\leq
k_{\alpha_0}a_{\alpha_0}\wedge v_x$, which implies that
$n_{a_{\alpha_0}}a_{\alpha_0}\leq v_x\leq x\leq a^{'}_{\alpha_0}$,
a contradiction. We obtain that
$k_{\alpha_0}a_{\alpha_0}\wedge(v_x\vee(x\ominus
v_x)=k_{\alpha_0}a_{\alpha_0}\wedge(x\ominus v_x)$, hence
$k_{\alpha_0}a_{\alpha_0}\leq x\ominus v_x$.

Now, let $\alpha_0\in\Lambda$ with
$k_{\alpha_0}=ord(a_{\alpha_0})$. Because
$n_{a_{\alpha_0}}a_{\alpha_0}\in S(E)$ and $v_x$ is the greatest
sharp element under $x$, we have $n_{a_{\alpha_0}}a_{\alpha_0}\leq
v_x\leq x$. It follows that $v_x$ and $x\ominus v_x$ are
compatible with  $n_{a_{\alpha_0}}a_{\alpha_0}$. Evidently
$n_{a_{\alpha_0}}a_{\alpha_0}\wedge (x\ominus v_x)=0$, because
$n_{a_{\alpha_0}}a_{\alpha_0}\leq v_x.$

Since we have proved that for every $\alpha\in\Lambda$ elements
$k_{\alpha}a_{\alpha}, v_x$ and $x\ominus v_x$ are pairwise
compatible. By Lemma 2.2, (iv) we have
$$v_x=\bigvee\{k_{\alpha}a_{\alpha}\wedge
v_x|\alpha\in\Lambda\}=\bigvee\{k_{\alpha}a_{\alpha}|k_{\alpha}=ord(a_{\alpha}),
\alpha\in\Lambda\}$$ and $$x\ominus v_x=\bigvee
\{k_{\alpha}a_{\alpha}\wedge (x\ominus
v_x)|\alpha\in\Lambda\}=\bigvee\{k_{\alpha}a_{\alpha}|k_{\alpha}\neq
ord(a_{\alpha}), \alpha\in\Lambda\}.$$ This proves the Theorem.

{\bf Theorem 4.2}. Let $(E; \oplus, 0, 1)$ be a sharply dominating
Archimedean atomic lattice effect algebra. Then to every
$(o)$-continuous state $\omega$ on $S(E)$ there exists a state
$\hat{\omega}$ on $E$ such that $\hat{\omega}|E=\omega$.

{\bf Proof.} Assume that a map $\omega: S(E)\rightarrow [0, 1]$ is
an $(o)$-continuous state on $S(E)$ and let us consider a map
$\hat{\omega}: E\rightarrow [0, 1]$ as follows:

(1). For every atom $a\in E$ let
$\hat{\omega}(a)=\frac{\omega(n_aa)}{n_a}$, where $n_a=ord(a)$.

(2). For every $x\in E, x\neq 0$ with basic decomposition
$x=v_x\oplus(\bigoplus\{k_{\alpha}a_{\alpha}|{\alpha}\in\Lambda\})$
set

$\hat{\omega}(x)=\omega(v_x)+\sup\{\sum_{\alpha\in\cal
K}k_{\alpha}\hat{\omega}(a_{\alpha})|{\cal K}\subseteq \Lambda$
finite $\}.$

Let us show that if $x=\bigoplus\{l_{\beta}b_{\beta}|\beta\in\cal
B\}$, where $\{b_{\beta}|\beta\in\cal B\}$ is an arbitrary set of
atoms of $E$ and positive integers $l_{\beta}$ satisfy
$l_{\beta}\leq n_{b_{\beta}}, \beta\in\cal B$ then
$\hat{\omega}(x)=\sup\{\sum_{\beta\in\cal
K}l_{\beta}\hat{\omega}(b_{\beta})|{\cal K}\subseteq \cal B$
finite$\}.$ Set ${\cal B}_1$=$\{\beta\in\cal B$$|l_{\beta}=$ ord
$(b_{\beta})\}$ and ${\cal B}_2$=$\{\beta\in\cal
B$$|l_{\beta}\neq$ ord $(b_{\beta})\}$. Evidently ${\cal
B}_1$$\cap{\cal B}_2=\emptyset$ and ${\cal B}_1$$\cup{\cal
B}_2=\cal B$. By Theorem 4.1 we have

$x=(\bigoplus\{l_{\beta}b_{\beta}|\beta\in {\cal
B}_1\})\oplus(\bigoplus\{l_{\beta}b_{\beta}|b\in {\cal B}_2\}),
\bigoplus\{l_{\beta}b_{\beta}|\beta\in {\cal B}_1\}=v_x$ and
$x=v_x\oplus(\bigoplus\{l_{\beta}b_{\beta}|\beta\in {\cal B}_2\})$
is a basic decomposition of $x$.

Assume $K\subseteq \cal B$ finite and set $F=K\cap {\cal B}_1$,
$G=K\cap {\cal B}_2$ and $x_K=\bigoplus_{\beta\in
K}l_{\beta}b_{\beta}$, $x_F=\bigoplus_{\beta\in
F}l_{\beta}b_{\beta}$ and $x_G=\bigoplus_{\beta\in
G}l_{\beta}b_{\beta}$. Then $x_K=x_F\oplus x_G$ is the basic
decomposition of $x_K$ by Theorem 4.1 and hence by (2) we have
$\hat{\omega}(x_K)=\omega(x_F)+\hat{\omega}(x_G)$. Further, by
Theorem 4.1 and the definition of $\bigoplus D$ for
$\oplus$-orthogonal $D$ we have $x_K\uparrow x, x_F\uparrow v_x$
and $x_G\uparrow x\ominus
v_x=\bigoplus\{l_{\beta}b_{\beta}|\beta\in \cal B$$_2\}$. Since
$\omega$ is $(o)$-continuous on $S(E)$, using the definition of
$\hat{\omega}$ we obtain that $\sup\{\hat{\omega}(x_K)|K\subseteq
\cal B$ finite $\}=\sup\{\omega(x_F)+\hat{\omega}(x_G)|F=K\cap
\cal B$$_1, G=K\cap \cal B$$_2, K\subseteq\cal B$ finite
$\}=\sup\{\omega(x_F)|F\subseteq {\cal B}$$_1$ finite
$\}+\sup\{\hat{\omega}(x_G)|G\subseteq {\cal B}_2$ finite
$\}=\omega(v_x)+\sup\{\sum_{\beta\in
G}l_{\beta}\hat{\omega}(b_{\beta})|\beta\in\cal
B$$_2\}=\hat{\omega}(x)$.

Now, let $x, y\in E$ with $x\leq y^{'}$. Then there exist sets
$\{a_{\alpha}|\alpha\in \Lambda\}$ and $\{b_{\beta}|\beta\in \cal
B\}$ of atoms of $E$ such that
$$x=\bigoplus\{k_{\alpha}a_{\alpha}|\alpha\in \Lambda\},
y=\bigoplus\{l_{\beta}b_{\beta}|\beta\in \cal B\}.$$ Further, as
was shown above, $\hat{\omega}(x)= \sup\{\sum_{\alpha\in
K}k_{\alpha}\hat{\omega}(a_{\alpha})|K\subseteq\Lambda$ finite$\}$
and $\hat{\omega}(y)= \sup\{\sum_{\beta\in
M}l_{\beta}\hat{\omega}(b_{\beta})|K\subseteq\cal B$ finite$\}$,
it follows that $\hat{\omega}(x\oplus
y)=\hat{\omega}(x)+\hat{\omega}(y)$, since $x\oplus
y=(\bigoplus\{k_{\alpha}a_{\alpha}|\alpha\in \Lambda\})\oplus
(\bigoplus\{l_{\beta}b_{\beta}|\beta\in \cal B\})$. Since $0, 1\in
S(E)$, we have $\hat{\omega}(0)=\omega(0)=0$ and
$\hat{\omega}(1)=\omega(1)=1$. This proves that $\hat{\omega}$ is
a state on $E$. Moreover, $\hat{\omega}|S(E)=\omega$ because $x\in
S(E)$ iff $x=v_x$.

{\bf Remark 4.3}. The next example shows that in Theorem 4.2 the
assumption that $E$ is lattice ordered cannot be omitted.

{\bf Example 4.4}. The smallest effect algebra $E$ admitting no
states has been presented in [14], Example 2.3. Namely, $E=\{0, a,
b, c, 2a, 2b, 2c, 3b, 1\}$ and $1=a\oplus b\oplus c=3a=4b=3c$,
which gives $b\oplus c=2a, a\oplus b=2c, a\oplus c=3b$. Thus, $E$
is not lattice ordered, because $a\vee b$ does not exists. Further
$S(E)=\{0, 1\}$, hence $E$ is sharply dominating and Archimedean.
Nevertheless, there is no states on $E$ extending a state $\omega$
existing on $S(E)$.

Finally, let us note that every complete (hence every finite)
lattice effect algebra $E$ is Archimedean (see [12], Theorem 3.3).
Moreover, $E$ is sharply dominating, because $S(E)$ is a complete
sublattice of $E$ and hence $\bigwedge Q\in E$ for all $Q\subseteq
S(E)$ (see [7], Theorem 3.7). Further, it is easy to verify that a
direct product of Archimedean sharply dominating atomic lattice
effect algebras is again an atomic, Archimedean, sharply
dominating and lattice ordered effect algebra.

\vskip 0.2 in

{\bf Acknowledgments}

\vskip 0.2 in

This project is supported by Natural Science Fund of China
(10471124), by Natural Science Fund of Zhejiang Provincial
(M103057) and by APVV Sino-Slovak Scientific and Technological
Cooperation, SK-CN-00906, and by VEGA 1/3025/06, M\v{S} SR.

\vskip 0.1 in

{\bf References}

\vskip 0.2 in

\noindent 1. C. C. Chang. Algebraic analysis of many-valued
logics. Trans. Amer. Math. Soc. 88(1958), 467-490

\noindent 2. D. J. Foulis, M. K. Bennett. Effect algebras and
unsharp quantum logics, Found. Phys. 24(1994), 1331-1352.

\noindent 3. R. J. Greechie. Orthomodular lattices admitting no
states. J. Combin. Theory A 10(1971), 119-132

\noindent 4. S. P. Gudder. $D$-algebras. Found. Phys. 26(1996),
813-822

\noindent 5. S. P. Gudder. Sharply dominating effect algebras.
Tatra Mt. Math. Publ. 15 (1998), 23-30

\noindent 6. S. P. Gudder. S-dominating effect algebras. Internat.
J. Theoret. Phys. 37 (1998), 915-923

\noindent 7. G. Jen\v{c}a, Z. Rie\v{c}anov\'a. On sharp elements
in lattice ordered effect algebras. BUSFFAL. 80(1999), 24-29

\noindent 8. F. K\^{o}pka, F. Chovanec. D-posets. Math. Slovaca,
43(1994), 21-34

\noindent 9. F. K\^{o}pka. Compatibility in $D$-posets. Internat.
J. Theor. Phys. 34(1995), 1525-1531

\noindent 10. Z. Rie\v{c}anov\'a. On proper orthoalgebras,
difference posets and Abelian relative inverse semigroups. Tatra
Mountains Math. Publ. 10(1997), 119-128

\noindent 11. Z. Rie\v{c}anov\'a. Macneille completions of
$D$-posets and effect algebras. Internat. J. Theoret. Phys.
39(2000), 855-865

\noindent 12. Z. Rie\v{c}anov\'a. Archimedean and block-finite
lattice effect algebras. Demonstratio Math. 33(2000), 443-452

\noindent 13. Z. Rie\v{c}anov\'a. Generalization of blocks for
$D$-lattice and lattice-ordered effect algebras. Internat. J.
Theoret. Phys. 39(2000), 231-236

\noindent 14. Z. Rie\v{c}anov\'a. Proper effect algebras admitting
no states. 40(2001), 1683-1691

\noindent 15. Z. Rie\v{c}anov\'a. Lattice effect algebras with
($o$)-continuous faithful valuations. Fuzzy Sets and System.
124(2001), 321-327

\noindent 16. Z. Rie\v{c}anov\'a. Smearings of states defined on
sharp elements onto effect algebras. Internat. J. Theoret. Phys.
41(2002), 1511-1524

\noindent 17. Z. Rie\v{c}anov\'a. Continuous lattice effect
algebras admitting order-continuous states. Fuzzy Sets and System.
136(2003), 41-54

\noindent 18. Z. Rie\v{c}anov\'a. Basic decomposition of elements
and Jauch-Piron effect algebras. Fuzzy Sets and System. 155(2005),
138-149

\end{document}